\documentclass[12pt,leqno]{article}
\usepackage{color}

\def\pf{\proof}

\def\epf{\hfill$\Box$ \medskip}

\usepackage{tikz}
\usepgflibrary{shapes.geometric}
\usepackage{amsmath,amsthm, amscd, amssymb, amsfonts, mathrsfs}
\usepackage[all]{xy}
\usepackage{color}
\usepackage[vcentermath]{youngtab}

\usepackage{graphics}
\usepackage{amssymb}
\usepackage{amscd}
\usepackage{amsmath}
\input xy 
\xyoption{all}
\title{Lusztig's strata are locally closed}
\newtheorem{theorem}{Theorem}[section]
\newtheorem{lemma}[theorem]{Lemma}

\author{Giovanna Carnovale\\
Dipartimento di Matematica ``Tullio Levi-Civita''\\
Torre Archimede - via Trieste 63 - 35121 Padova - Italy\\
email: carnoval@math.unipd.it }
\date{}

\begin{document}
\maketitle
\begin{abstract}
Let $G$ be a connected reductive algebraic group over an algebraically closed field $k$. We consider the strata in $G$ defined by Lusztig as fibers of a map given in terms truncated induction  of Springer representations.
We extend to arbitrary characteristic the following two results:  Lusztig's strata are locally closed and the irreducible components  of a stratum $X$ are those sheets for the $G$-action on itself that are contained in $X$.
\end{abstract}

\section{Introduction}
The present paper answers a question by G. Lusztig on the extension to arbitrary characteristic of results in \cite{gio-MR} concerning strata as defined in \cite[\S2]{lustrata}. 
Lusztig's strata in a connected reductive group $G$ are defined as fibers of a map from $G$ to the set of irreducible representations of its Weyl group, constructed in terms of truncated induction of Springer's representations for trivial local systems. It is  observed in  \cite[\S2]{lustrata} that they are a $G$-stable union of orbits of the same dimension and are unions of Jordan classes, i.e., the locally closed subsets introduced in \cite[\S3]{lusztig-inventiones} that provide the stratification with respect to which character sheaves are smooth. 
It does not immediately follow from the definition that these fibers have  topological or geometric properties.  However, it was proved in  \cite{gio-MR} when the characteristic of the base field is good, that strata are locally closed and that they are unions of sheets for the action of $G$ on itself by conjugation. In this paper we show that these results holds in arbitrary characteristics.  The proof of all both statements relies on the observation that if a Jordan class lies in a stratum, then also the regular part of its closures lies in the stratum. As a consequence, we conclude that the irreducible components of a stratum are precisely the sheets contained therein. This implies that whenever two sheets have non-empty intersection, then the stratum containing both is singular. This happens, for instance, when the root system of $G$ is not simply-laced and the stratum contains the subregular unipotent conjugacy class. 

\section{Notation}Let $G$ be a connected reductive algebraic group over an algebraically closed field $k$. We denote by $G\cdot x$ the $G$-conjugacy class of $x\in G$, whereas $G^{x\circ}$ denotes the identity component of the stabiliser of $x$. 
For $m\in{\mathbb N}$ we set $G_{(m)}=\{x\in G~|~\dim G\cdot x=m\}$. 
Since the dimension of the stabiliser of a point is an upper semicontinuous function \cite[\S2, page 7]{GIT} we have $\bigcup_{m\leq d}\overline{G_{(m)}}=\bigcup_{m\leq d}G_{(m)}$ so 
$G_{(>d)}:=\bigcup_{m\geq d+1}G_{(m)}$ is open in $G$.\footnote{Observe that in \cite[page 8]{gio-espo} it is erroneously stated that $\overline{G_{(m)}}\bigcup_{m\leq d}G_{(m)}$: this does not affect any of the following statements in that paper.} 
 
For a set $Y\subset G$ we define $d_Y:=\max\{d\in{\mathbb N}~|~G_{(d)}\cap Y\neq\emptyset\}$ and  $Y^{reg}:=Y\cap G_{(d_Y)}$. 
The sets $G_{(m)}=(\bigcup_{l\leq m}G_{l})\cap G_{(>m-1)}$ are locally-closed and their irreducible components are called the {\em sheets} of the $G$-action. 
We will denote by $W$ the Weyl group of $G$, by ${\rm Irr}(W)$ the set of isomorphism classes of the complex irreducible representations of $W$. For $s$ in a maximal torus $T$ of $G$, we set $W_s=N_G(T)\cap G^{s\circ}/T$, which is the Weyl group of $G^{s\circ}$. When we write $g=su$ we mean that $su$ is the Jordan decomposition of $g$ with $s$ semisimple and $u$ unipotent. 

\subsection{Jordan classes and sheets}We recall from  \cite[\S3]{lusztig-inventiones} that the group $G$ is the disjoint union of finitely many, locally closed, smooth, irreducible, $G$-stable sets, each contained in some $G_{(m)}$,
 which we call Jordan classes and can be described as follows: the Jordan class containing $g=su$ is 
$J(su)=G\cdot((Z(G^{s\circ})^\circ s)^{reg}u)$. In other words, an element with Jordan decomposition $rv$ lies in $J(su)$ if and only if it is conjugate to some $s'u$ with $G^{s\circ}=G^{s'\circ}$ and $s'\in Z(G^{s\circ})^\circ s$.  Last condition guarantees that $J(su)$ is irreducible.  

Clearly, if $J(su)\subset G_{(d)}$, then $\overline{J(su)}\subset\bigcup_{d\leq m}G_{(d)}$ and 
$J(su)\subset\overline{J(su)}^{reg}=\overline{J(su)}\cap G_{(d)}$. 
Next Lemma is a combination of  \cite[Propositions 4.5 and 4.7]{gio-espo}: both proofs hold in any characteristic and can be repeated verbatim.
\begin{lemma}Let $J(su)\subset G_{(d)}$ be the Jordan class of $g=su$. Then
\begin{equation*}\overline{J(su)}=\bigcup_{z\in Z(G^{s\circ})^\circ}\overline{G\cdot sz {\rm Ind}_{G^{s\circ}}^{G^{zs\circ}}(G^{s\circ}\cdot u)}\subseteq \bigcup_{l\leq d}G_{(l)}
\end{equation*}
and
\begin{equation*}\overline{J(su)}^{reg}=\overline{J(su)}\cap G_{(d)}=\bigcup_{z\in Z(G^{s\circ})^\circ}G\cdot sz {\rm Ind}_{G^{s\circ}}^{G^{zs\circ}}(G^{s\circ}\cdot u).\end{equation*}
\end{lemma}
Let $S$ be a sheet contained in  $G_{(d)}$. Since the Jordan classes contained in $G_{(d)}$ are irreducible they are all contained in a sheet, and since they are finitely-many, there is a Jordan class $J$ such that $\overline{J}=\overline{S}$, so $S=\overline{S}^{reg}=\overline{J}^{reg}$.

\subsection{Lusztig's strata}We recall from \cite[\S2]{lustrata} the map $\phi_G\colon G\to {\rm Irr}(W)$ defined on $g=su$ as:
$\phi_G(g)={\bf j}_{W_s}^W\rho_u^{W_s}$, where $\rho_u^{W_s}$ is Springer's respresentation of $W_s$ associated with $u$ and trivial local system and ${\bf j}_{W_s}^W$ is Lusztig-Spaltenstein's induction \cite{lusp}.
The fibers of this map are $G$-stable, each of them is contained in some  $G_{(d)}$ and by construction they are union of Jordan classes.

\begin{theorem}Let $X$ be a stratum in $G$. Then 
\begin{enumerate}
\item[(1)]\label{item:uno}$X$ is locally closed.
\item[(2)]\label{item:due} $X$ is a union of sheets.
\item[(3)]\label{item:tre} The sheets contained in $X$ are its irreducible components. 
\end{enumerate}
\end{theorem}
\pf 1. We claim  that if $J(su)\subset X$, then $\overline{J(su)}^{reg}\subset X$.
Let $rv\in \overline{J(su)}^{reg}$. Without loss of generality we may assume $r=zs$  for $z\in Z(G^{s\circ})^\circ$ and $v\in {\rm Ind}_{G^{s\circ}}^{G^{sz\circ}}(G^{s\circ}\cdot u)$.  
Then, $\phi_G(rv)={\bf j}_{W_{zs}}^W\rho^{W_{zs}}_v$. By \cite[Theorem 3.5]{lusp}, see also \cite[\S 4.1]{spalt} we have 
\begin{equation*}\phi_G(rv)=\phi_G(zsv)={\bf j}_{W_{zs}}^W\rho_v^{W_{zs}}={\bf j}_{W_{zs}}^W{\bf j}_{W_{s}}^{W_{zs}}\rho_u^{W_{s}}={\bf j}_{W_s}^W\rho_u^{W_{s}}=\phi_G(su).\end{equation*}
Therefore $X$ is a union of finitely many sets of the form $\overline{J(s_lu_l)}^{reg}$. If $X\subset G_{(d)}$, then $\overline{J(s_lu_l)}^{reg}\subset G_{(d)}$ for every $l$ so
\begin{equation*}X=\bigcup_l\overline{J(s_lu_l)}^{reg}=\overline{\bigcup_l J(s_lu_l)}^{reg}=\overline{\bigcup_l J(s_lu_l)}\cap G_{(d)}=
\overline{\bigcup_l J(s_lu_l)}\cap G_{(>d)}\end{equation*}
with $\overline{\bigcup_l J(s_lu_l)}$ closed and $G_{(>d)}$ open.
2 and 3. For any $d$, inclusion induces a partial ordering on the collection of the finitely-many irreducible sets of the form $\overline{J(su)}^{reg}$ contained in $G_{(d)}$. The maximal elements are the irreducible components of  $G_{(d)}$, i.e., the sheets of $G$ corresponding to the dimension $d$. Hence $X=X\cap G_{(d)}$ is  a union of sheets and these are the maximal irreducible subsets contained in $X$, i.e., its irreducible components.
\epf

\section{Acknowledgements}
We thank Prof. G. Lusztig for asking the question leading to this short note and for clarifying a doubt on compatibility of truncated induction. We also thank F. Esposito for pointing out reference \cite{GIT} and for a discussion on the colsures of the sets $G_{(d)}$.  The author acknowledges support from by DOR1898721/18,  DOR1717189/17  and BIRD179758/17 funded by the University of Padova.

\end{document}